\documentclass[12pt,oneside,reqno]{amsart}
\usepackage[usenames]{color}
\usepackage[font=small,labelfont=bf]{caption}
\usepackage{sidecap}
\usepackage{float}
\usepackage{upgreek}
\usepackage{mathabx}
\usepackage{mathptmx}
\usepackage{bookman}

\usepackage{hyperref}

\usepackage{amssymb}
\usepackage{amscd}
\usepackage{array}
\usepackage{verbatim}

\DeclareSymbolFont{rmlargesymbols}{OMX}{mdbch}{m}{n}
\DeclareMathSymbol{\rmintop}{\mathop}{rmlargesymbols}{82}
\DeclareMathSymbol{\rmointop}{\mathop}{rmlargesymbols}{72}
\DeclareMathSymbol{\rmsumop}{\mathop}{rmlargesymbols}{80}
\DeclareMathSymbol{\rmunionop}{\mathop}{rmlargesymbols}{83}
\DeclareMathSymbol{\rmintersectop}{\mathop}{rmlargesymbols}{84}
\DeclareMathSymbol{\rmtensorop}{\mathop}{rmlargesymbols}{79}
\DeclareMathSymbol{\rmdirectsumop}{\mathop}{rmlargesymbols}{77}

\def\ccite#1{\textcolor{red}{\cite{#1}}}
\numberwithin{equation}{section}

\definecolor{mycolor}{rgb}{0.39,0.26,0.13}
\definecolor{mygray}{rgb}{0.1,0.1,0.1}
\definecolor{myblue}{rgb}{0,0,0.2}
\definecolor{whiteblue}{rgb}{0.8,1,0.9}

\begin{document}


\title[A P\lowercase{roof} O\lowercase{f} C\lowercase{antor's} 
T\lowercase{heorem}]
{\Large\rm A P\lowercase{roof} O\lowercase{f} C\lowercase{antor's} 
T\lowercase{heorem}}

\author{S. W\lowercase{alters}}
\date{March 30, 2020}
\address{Department of Mathematics \& Statistics, University  of Northern B.C., Prince George, B.C. V2N 4Z9, Canada.}
\email[]{samuel.walters@unbc.ca}
\subjclass[2000]{42A05; 42A16, 11L03, 33B10, 42A55} 
\keywords{Trigonometric series, trigonometric polynomials, Cantor/Lebesgue theorem, Fourier series}
\urladdr{http://hilbert.unbc.ca}

\begin{abstract}
We present a short proof of Cantor's Theorem (circa 1870s): if $a_n \cos nx + b_n \sin nx \to 0$ for each $x$ in some (nonempty) open interval, where $a_n, b_n$ are sequences of complex numbers, then $a_n$ and $b_n$ converge to 0.
\end{abstract}

\maketitle

\section{Proof Of Cantor's Theorem}

\noindent{\bf Cantor's Theorem.} Let $a_n, \, b_n$ be sequences of complex numbers such that
\[
\lim_{n\to\infty} \ a_n \cos nx + b_n \sin nx \ = \ 0
\]
for each $x$ in some open interval $(c,d)$. Then $a_n \to 0$ and $b_n \to 0$.
\bigskip

The proof presented here consists of reduction to the case $C_n  \sin nx \to 0,$ which is covered by Lemma B below and which we proceed to prove first.

\medskip

\noindent{\bf Lemma A.} Let $\delta > 0 $ be given and let $B_n$ be a {\it bounded} sequence of complex numbers such that $B_n \sin n t \to 0$ as $n\to \infty$ for each $0 < t < \delta$. Then $B_n \to 0$. 

\begin{proof} Replacing $t\in(0,\delta)$ by $2t$ we have 
\[
B_n \sin (n 2t) = 2 B_n \sin(nt) \cos(nt) \to 0
\] 
so the limit $B_n \sin n x \to 0$ holds also for $0 < x < 2\delta$. By repeating this doubling procedure a finite number of times, we get $B_n \sin n x \to 0$ for $0 \le x \le 2\pi \le N\delta$ for some fixed positive integer $N$. Therefore, $f_n(x) := B_n \sin n x \to 0$ for each $x$ in the closed interval $[0,2\pi]$. Further, since the sequence $B_n$ is bounded, the functions $f_n$ are bounded by an (integrable) constant on $[0,2\pi]$. Therefore, by the Dominated Convergence Theorem (\ccite{Rudin}, Theorem 11.32)
\[
0 \ = \ \lim_n \rmintop_0^{2\pi} |f_n(x)| dx \ =\  \lim_n \ |B_n|  \rmintop_0^{2\pi} |\sin nx| dx 
\ = \  \lim_n \ 4 |B_n| 
\]
since $\int_0^{2\pi} |\sin nx| dx = 4$ for each integer $n\ge1$. It follows that $B_n \to 0$.
\end{proof}

The next lemma removes the boundedness condition in Lemma A.
\medskip

\noindent{\bf Lemma B.} Let $\delta > 0 $ be given and let $C_n$ be a sequence of complex numbers such that $C_n \sin n t \to 0$ as $n\to \infty$ for each $0 < t < \delta$. Then $C_n \to 0$. 
\begin{proof} 
Consider the bounded sequence $B_n = \frac{ |C_n| }{ 1+|C_n|},$ for which we have
\[
\lim_{n} \ B_n \sin nx \ =\ \lim_{n} \ \frac{ |C_n| \sin nx}{ 1+|C_n|}  = 0
\]
holds for each $x$ in $(0,\delta)$. Hence $B_n \to 0$ by Lemma A, which in particular means that there is an integer $M$ such that $B_n < \tfrac12$ for $n > M$. This says that $|C_n| < 1$ (for $n>M$) so that $C_n$ is bounded. By Lemma A, applied again to $C_n,$ we see that $C_n\to 0$. 
\end{proof}
\medskip

\noindent{\it Proof of Cantor's Theorem.} We are given that
\[
a_n \cos nx + b_n \sin nx \to 0
\]
for each $c < x < d$. Let $\delta > 0$ such that $x + \delta < d,$ so that 
$c < x+t < d$ holds for each $0 \le t < \delta$. Then since the zero limit holds also for $x + t,$ we have
\begin{align}\label{lastline}
a_n \cos n&(x+t) + b_n \sin n(x+t) = a_n \cos(nx+nt) + b_n \sin(nx+nt)	\notag
\\
&= a_n \big[ \cos(nx) \cos(nt) - \sin(nx)\sin(nt)  \big] + b_n \big[ \sin(nx) \cos(nt) + \cos(nx) \sin(nt)  \big]	\notag
\\
&= \big[ a_n \cos(nx) + b_n \sin(nx) \big] \cos(nt) 
+ \big[ b_n \cos(nx) - a_n \sin(nx) \big] \sin(nt).
\end{align}
Since this whole sum goes to 0, and since the first term in \eqref{lastline} goes to 0 by hypothesis, it follows that the second term goes to 0, i.e. $C_n \sin(nt) \to 0$ for each $0 < t < \delta$ where $C_n := b_n \cos(nx) - a_n \sin(nx)$ is a fixed sequence (independent of $t$) since here $x$ is fixed. Lemma B above now applies to this and gives
\[
C_n = b_n \cos(nx) - a_n \sin(nx) \ \to \ 0
\]
as $n\to\infty$. Therefore taking the sum of absolute squares we get
\[
\big|a_n \cos(nx) + b_n \sin(nx)\big|^2 + 
\big|b_n \cos(nx) - a_n \sin(nx)\big|^2 = |a_n|^2 + |b_n|^2 
\] 
which has zero limit, hence the result. \qed

\medskip


\end{document}